\documentclass[12pt]{article}
\usepackage{epsfig}
\usepackage{amsmath}
\usepackage{amsfonts}
\usepackage{amscd}

\begin{document}

\title{A cohomological description of Abelian bundles and gerbes
 \footnote{Submitted to the proceedings of the
 XXth Workshop on Geometric Methods in Physics, Bialowie\.{z}a, July
 1-7, 2001,}
}
\author{
{ Roger Picken}\\ Departamento de Matem\'{a}tica
and\\
CEMAT,
Centro de Matem\'{a}tica e Aplica\c{c}\~{o}es\\
Instituto
Superior T\'{e}cnico\\ Av. Rovisco Pais\\ 1049-001 Lisboa\\ Portugal\\
e-mail: rpicken@math.ist.utl.pt}
\date{April 9, 2003}
\maketitle
\begin{abstract}
We describe the geometrical ladder of equations for Abelian bundles and
gerbes, as well as higher generalisations, in terms of the cohomology of an
operator that combines de Rham and \v{C}ech cohomology.

\end{abstract}
\section{Introduction}
Gerbes with connection appear in differential geometry as a natural
higher-order generalization of abelian bundles with connection, and thus,
from the physics standpoint, provide a possible framework in which to
generalise Abelian gauge theory. They first appeared in algebraic geometry
\cite{Gir}, and were subsequently developed by Brylinski
\cite{Bry}, whose motivation was to generalise the geometrical
interpretation of the second integral cohomology of a manifold $M$, $H^2(M,
{\bf Z})$, in terms of the curvature of a complex line bundle, to $H^3(M,
{\bf Z})$. There has been renewed interest in the subject recently
following a concrete approach due to Hitchin and Chatterjee~\cite{Hi99},
and due to the appearance of possible applications in physics, for instance in anomalies \cite{cmm}, new geometrical structures in string theory
\cite{wit} and Chern-Simons theory \cite{gom}.

Bundles and gerbes, as well as higher generalisations ($n$-gerbes), can be
understood both in terms of local geometry, i.e. local functions and forms,
and in terms of non-local geometry, i.e. holonomies and parallel
transports, and these two viewpoints are equivalent, in a sense made
precise by Mackaay and the author in \cite{MP}, following on from work by
Barrett \cite{Bar91} and Caetano and the author \cite{CP94}. For gerbes,
holonomy and parallel transport are along embedded surfaces in the
manifold, instead of along loops or paths (as indeed might be expected, since for
gerbes everything is ``one dimension up'', compared to bundles). The
non-local viewpoint was further explored by the author in \cite{P03}, as a
special case of a topological quantum field theory framework introduced in
\cite{PS}. Here we wish to examine in greater detail the geometrical
ladders of local functions and forms that appear in the local geometry
perspective. We use a simple cohomological approach that is similar to Deligne
cohomology.

The article is organized as follows. In section 2 we recall the equations
governing principal $U(1)$ bundles with connection in terms of transition
functions and local connection $1$-forms.
In section 3 we generalise to $U(1)$ gerbes with connection and give a
simple example of a gerbe on the $3$-sphere.
In
section 4 we present our cohomological framework for describing the
whole ladder of $n$-gerbes and their equivalences.
Section 5 contains some comments.

\vskip0.5cm
\centerline{\bf Acknowledgements}
\vskip.5cm
\noindent
I am very grateful to Marco Mackaay for much stimulating collaboration, and
for suggesting the method of dealing with the signs in section
\ref{cohom}.
It is a pleasure to thank the organizers of the XXth Bialowie\.{z}a
Workshop for giving me the opportunity to present this material for
the first time, and for their remarkable achievement in creating and
maintaining the excellent scientific and social traditions of the
Bialowie\.{z}a meetings.

This work was supported by {\em Programa Operacional
``Ci\^{e}ncia, Tecnologia, Inova\c{c}\~{a}o''} (POCTI) of the
{\em Funda\c{c}\~{a}o para a Ci\^{e}ncia e a Tecnologia} (FCT),
cofinanced by the European Community fund FEDER.

\label{Intro}
\section{The local equations for bundles with connection}
We start by recalling a few well-known facts about principal bundles. A principal $G$-bundle $P$ over a manifold $M$ is given by a projection map $\pi:P\rightarrow M$, where $\pi^{-1}(x)$ is called the fibre over $x\in M$, $P$ is called the total space, and $M$ is called the base space, together with a (right, effective) $G$-action on $P$, written $p.g$ for $p\in P$ and $g\in G$, preserving fibres: $\pi(p.g)=\pi(g)$. Here $G$ is a Lie group, and all manifolds and maps are smooth. The axiom of local triviality says that for any $x\in M$ there exists an open neighbourhood $U$ of $x$ such that $\pi^{-1}(U)$ is isomorphic to $U\times G$. A local trivialization is given by a local section, i.e. $s:U\rightarrow P$, satisfying $\pi(s(x))=x$ for all $x\in U$, via $p=s(x).g \in \pi^{-1}(U)\leftrightarrow (x,g)\in U\times G$.

Given two open sets $U_1$ and $U_2$ which intersect, and local sections $s_1$ and $s_2$ on $U_1$ and $U_2$ respectively, they are related on the overlap $U_1\cap U_2$ by a transition function $g_{12}$ defined by:
\[
s_2(x)=s_1(x). g_{12}(x), \quad \forall x\in U_1\cap U_2.
\]
If we introduce a third local section $s_3$, defined on an open set $U_3$ that intersects $U_1\cap U_2$ non-trivially, then, by writing $s_3(x)$ in two different ways, we obtain the equation
\[
g_{12}(x) g_{23}(x) = g_{13}(x), \quad \forall x\in U_1\cap U_2\cap U_3
\]
on the triple overlap. In particular we have, choosing $U_3=U_1$,
\begin{equation}
g_{21}(x) = g_{12}^{-1}(x).
\label{ginverse}
\end{equation}

The (Ehresmann) connection is a $1$-form $\omega$ on $P$ taking values in the Lie algebra of $G$, and satisfying some extra properties. Using local sections we obtain $1$-forms on $M$ by pull-back. Let $A_1=s_1^\ast \omega$ and
$A_2=s_2^\ast \omega$ be two such local $1$-forms, defined on $U_1$ and $U_2$ which intersect non-trivially. Then they are related by the equation
\[
A_2=g_{12}^{-1}A_1 g_{12} + g_{12}^{-1}d g_{12}
\]
on $U_1\cap U_2$, due to the properties of $\omega$.

If we use different local sections ${s_i}':U_i\rightarrow P$ instead, related to the original local sections $s_i$ by
\[
{s_i}'(x) = s_i(x).h_i(x), \quad \forall x\in U_i,
\]
where $h_i:U_i\rightarrow G$, we get gauge equivalent transition functions and connection $1$-forms:
\begin{eqnarray*}
g_{ij}' & = & h_i^{-1} g_{ij} h_j, \\
A_i' & = & h_i^{-1} A_i h_i + h_i^{-1} d h_i.
\end{eqnarray*}

We now specialize to the case when $G$ is $U(1)$, i.e. Abelian. We will replace functions with values in $G$ by their logarithms, so that all equations take values in $i\bf R$, the Lie algebra of $G$, and are modulo integer multiples of $2\pi i$. We will assume that all the open sets, and multiple overlaps $U_{ijk\dots}:= U_i\cap U_j \cap U_j \dots $ in our cover of $M$ are contractible. A $G$-bundle with connection is then given by transition functions
\[
\ln g_{ij}: U_{ij}\rightarrow i{\bf R},
\]
antisymmetric under exchange of indices because of equation~(\ref{ginverse}), and connection $1$-forms
\[
A_i\in \Lambda^1(U_i),
\]
satisfying
\begin{equation}
\ln g_{jk} - \ln g_{ik} + \ln g_{ij} =0
\label{b1}
\end{equation}
on $U_{ijk}$ (we will see the reason for writing the terms in this order in section~\ref{cohom}) and
\begin{equation}
i(A_j-A_i) =d\ln g_{ij}
\label{b2}
\end{equation}
on $U_{ij}$. We may also introduce the curvature $2$-form $F$ on $M$ defined by
\begin{equation}
F=dA_i
\label{b3}
\end{equation}
on each open set $U_i$ (which is indeed globally defined because of equation~(\ref{b2})). Finally $F$ satisfies the Bianchi identity
\begin{equation}
dF=0
\label{b4}
\end{equation}
on $M$.
Transition functions $\ln g'_{ij}$ and connection $1$-forms $A'_i$ are gauge-equivalent to transition functions $\ln g_{ij}$ and connection $1$-forms $A_i$, iff there exist functions
\[
\ln h_i: U_{i}\rightarrow i{\bf R},
\]
satisfying
\begin{equation}
\ln g'_{ij} - \ln g_{ij}  =  \ln h_j - \ln h_i
\label{be1}
\end{equation}
on $U_{ij}$ and
\begin{equation}
i(A'_i-A_i) = d\ln h_i
\label{be2}
\end{equation}
on $U_{i}$.

In this final formulation of the equations describing $U(1)$-bundles, there is no reference any longer to the total space $P$ of the principal bundle---the whole description is ``downstairs''. This is advantageous for the generalisation to gerbes with connection in the next section, where the notion of a total space is lacking.

\label{locbundle}
\section{The generalisation to gerbes with connection and an example}
The main guiding principle for understanding the generalisation from bundles to gerbes, is that for gerbes everything is one step up compared to bundles, in the form degree, the number of open sets in an overlap, or the dimension. Thus gerbes have transition functions defined on triple overlaps, a curvature $3$-form, and parallel transport defined along surfaces inside $M$, instead of paths. A good place to read more about the general background and geometrical applications of gerbes is in Hitchin's lectures \cite{Hi99}.

The data and equations defining a $U(1)$-gerbe with connection are analogous to the bundle case, except that there are now two separate layers of connections, connection $1$-forms and connection $2$-forms. A $U(1)$-gerbe with connection is given by
transition functions
\[
\ln g_{ijk}: U_{ijk}\rightarrow i{\bf R},
\]
completely antisymmetric under exchange of indices, connection $1$-forms
\[
A_{ij}\in \Lambda^1(U_{ij}),
\]
antisymmetric under exchange of indices, and connection $2$-forms
\[
F_{i}\in \Lambda^2(U_{i}),
\]
satisfying
\begin{equation}
\ln g_{jkl} - \ln g_{ikl} + \ln g_{ijl} - \ln g_{ijk} =0
\label{g1}
\end{equation}
on $U_{ijkl}$,
\begin{equation}
i(A_{jk} - A_{ik} + A_{ij}) = - d\ln g_{ijk}
\label{g2}
\end{equation}
on $U_{ijk}$, and
\begin{equation}
F_j - F_i = d A_{ij}
\label{g3}
\end{equation}
on $U_{ij}$. Again all equations are taken modulo $ 2\pi i$.
As before we may introduce a curvature form on $M$, this time a curvature $3$-form $G$, defined by
\begin{equation}
G=dF_i
\label{g4}
\end{equation}
on each open set $U_i$ (which is again globally defined because of equation~(\ref{g3})). Finally $G$ satisfies the ``Bianchi'' identity
\begin{equation}
dG=0
\label{g5}
\end{equation}
on $M$.

The notion of gauge equivalence for gerbes is also a higher notion, since equivalences are specified not just by giving functions, but also $1$-forms. More precisely, transition functions $\ln g'_{ijk}$, connection $1$-forms $A'_{ij}$ and connection $2$-forms $F'_i$ are gauge-equivalent to
transition functions $\ln g_{ijk}$, connection $1$-forms $A_{ij}$ and connection $2$-forms $F_i$
iff there exist functions
\[
\ln h_{ij}: U_{ij}\rightarrow i{\bf R}
\]
and $1$-forms
\[
B_i \in \Lambda^1(U_{i})
\]
satisfying
\begin{equation}
\ln g'_{ijk} - \ln g_{ijk}  =  \ln h_{jk} - \ln h_{ik} + \ln h_{ij}
\label{ge1}
\end{equation}
on $U_{ijk}$,
\begin{equation}
i(A'_{ij}-A_{ij}) = - d\ln h_{ij} + i(B_j - B_i)
\label{ge2}
\end{equation}
on $U_{ij}$, and
\begin{equation}
i(F'_i-F_i) = dB_i
\label{ge3}
\end{equation}
on $U_{i}$.

There is clearly a pattern in the above equations, which will be elucidated in the next section. One aspect which is sometime found to be puzzling, is that the connection $1$-form $A_{ij}$ for gerbes is not defined everywhere on $M$, but only on the double overlaps of the cover. It may seem that $A_{ii}$ is defined on $U_i$, but this is identically zero because of the antisymmetry condition on the indices. However one should really think of $A_{ij}$ as being merely a subsidiary ``transition'' connection for the genuine gerbe connection $2$-form $F_i$, which is defined on every patch of $M$. 

We will conclude this section with a simple example of a gerbe with connection on the $3$-sphere, which is a natural generalisation of the familiar monopole bundle on the $2$-sphere\footnote{Meinrenken \cite{mei} has recently obtained gerbes with connection on compact simple Lie groups $G$, which are $G$-equivariant. The purpose of our example is different, however, so we do not consider $S^3$ as a Lie group.}.  In fact, let us start by considering the monopole bundle on $S^2$, covered by two patches $U_1$ and $U_2$, which intersect in a (not-too-wide) strip around the equator, isomorphic to $S^1\times ]0,1[$. (This is not contractible, but we could introduce extra patches to get contractible overlaps if necessary.) The transition function $g_{12}:U_{12}\rightarrow U(1)$ is chosen to be a winding number $1$ map from $S^1$ to $U(1)$, constant along the transversal direction (which is why we do not want the strip to be too wide). On $U_1$ we choose the connection $1$-form $A_1=0$, and on $U_2$ we choose $i A_2=\overline{ d \ln g_{12}}$, meaning it is equal to $d \ln g_{12}$ on the overlap $U_{12}$ and is continued in some manner to the rest of $U_2$, which is possible since $U_2$ is contractible. Then we have the desired equation for a bundle with connection (\ref{b2}):
\[
i(A_2-A_1) = d \ln g_{12}
\]
on $U_{12}$. The curvature $F$ of the monopole connection, given by equation (\ref{b3}), has support contained in $U_2$, and integrating $F/2\pi$ over $S^2$ gives $1$. 

Now we take the $3$-sphere $S^3$, and cover it with three patches $U_1$, $U_2$ and $U_3$. $U_3$ covers the equator and one half of $S^3$, and $U_1$ and $U_2$ together cover the equator and the other half of $S^3$, in such a way that $U_{13}$ and $U_{23}$ are isomorphic to $U_1\times ]0,1[$ and $U_2\times ]0,1[$ respectively, where $U_i$ refer to the patches from the monopole bundle case. The intersection between $U_3$ and the union of $U_1$ and $U_2$ is a (not-too-wide) spherical shell isomorphic to $S^2\times ]0,1[$. The intersection of all three open sets is isomorphic to $S^1\times ]0,1[^2$. We take the transition function $g_{132}$ on this triple overlap to be the winding number $1$ map $g_{12}$ from $S^1$ to $U(1)$ used previously, taken to be constant along the transversal directions. The $1$-form connection is given by: $A_{13}=A_1$ (the $1$-form for the monopole, constant along the direction transversal to $U_1$), $A_{23}=A_2$ (constant in the direction transversal to $U_2$) and $A_{12}=0$. We thus have equation (\ref{g2})
\[
i(A_{32} - A_{12} + A_{13}) = - d\ln g_{132}.
\]
Now we choose the $2$-form connection $F_i$ as follows: $F_1=F_2=0$ and $F_3= \overline{F}$, meaning it is equal to the monopole curvature $F$ on the overlaps $U_{13}$ and $U_{23}$, and is continued in some manner to the rest of $U_3$, which is possible since $U_3$ is contractible. We have thus also satisfied equation (\ref{g3})
\[
F_2 - F_1 = d A_{12}=0,  \quad F_3 - F_1 = d A_{13}=dA_1, 
\quad F_3 - F_2 = d A_{23}=dA_2.
\]
The curvature of the gerbe connection $G$, given by equation (\ref{g4}), has support contained in $U_3$, and integrating $G/2\pi$ over $S^3$ gives the integral of $F/2\pi$ over $S^2$, i.e. $1$, after applying Stokes' theorem. 

The ``gerbopole'' we have just described is therefore a natural, higher generalisation of the monopole, and its winding number, or charge, also derives from the winding number $1$ map from the circle to the group. In the next section we will see that this map itself also has an interpretation in the dimensional ladder containing bundles and gerbes with connection. 

\label{locgerbe}
\section{A cohomological formulation}

In this section we wish to unify the equations for bundles and gerbes with connection from the previous sections, and generate the whole dimensional ladder of $n$-gerbes, by using a cohomological formulation. Our approach can be viewed as a variant of Deligne cohomology, in that it blends together de Rham and \v{C}ech cohomology. Let us first define the cochain groups $\Lambda^{p,n} $, whose elements are collections of $p$-forms, valued in $i \bf R$, defined on each $n$-fold overlap of open sets of our fixed cover of $M$.  When $n=1$, a cochain consists of a $p$-form on each open set of the cover, and when $n=0$, a cochain is a single, globally-defined $p$-form on $M$. We write cochains in the form 
$C_{ijk\dots} $ where $ijk\dots$ ranges over the $n$-fold overlaps, or $C$ (no index) when $n=0$.

There are two natural operators acting on these cochain groups, namely the exterior derivative
\[
d: \Lambda^{p,n} \rightarrow \Lambda^{p+1,n}, \quad C_{ijk\dots} \mapsto d C_{ijk\dots} 
\]
and the \v{C}ech coboundary operator
\[
\delta: \Lambda^{p,n} \rightarrow \Lambda^{p,n+1}, \quad \delta C_{i_1\dots i_{n+1}}=
\sum_{\alpha = 1}^{n+1} (-1)^{\alpha + 1} C_{i_1\dots \hat{i_\alpha}\dots  i_{n+1}}.
\]
Both of these operators are nilpotent, but in order to combine them into a new nilpotent operator, they should anticommute, which is not the case. This can be remedied either by multiplying $\delta$ by $(-1)^p$, or by multiplying $d$ by $(-1)^n$. We choose the latter solution, and define
\[
\bar{d}: \Lambda^{p,n} \rightarrow \Lambda^{p+1,n}, \quad \bar{d} =  (-1)^n \, d,
\]
which satisfies $\delta \bar{d} + \bar{d} \delta =0$. 

Now we can define the cohomology that is relevant for describing bundles and gerbes. Let the cochain groups $\Lambda^{(k)}$ be given by
\[
\Lambda^{(k)}= \bigoplus_{p+n=k} \Lambda^{p,n},
\]
and the operator $D: \Lambda ^{(k)}\rightarrow \Lambda ^{(k+1)}$ be defined by:
\[
D=\delta - \bar{d},
\]
satisfying $D^2=0$. A bundle with connection may now be defined to be a $2$-cocycle $\cal B$:
\[
{\cal B}= \ln g_{ij} + iA_i -iF \in \Lambda ^{(2)}, \quad  D {\cal B}=0,
\]
which is equivalent to equations (\ref{b1}) to (\ref{b4}), and a gerbe with connection may be defined to be a $3$-cocycle $\cal G$
\[
{\cal G} = \ln g_{ijk} +iA_{ij} +iF_i - iG \in \Lambda ^{(3)}, \quad D {\cal G}=0,
\]
which is equivalent to equations (\ref{g1}) to (\ref{g5}), as may be easily verified. 

Having simplified the equations in this way, we can extend the definition to cocycles of any order. Thus we define an $n$-gerbe with connection to be:
\[
{\cal H} \in \Lambda ^{(n+2)}, \quad D{\cal H}=0,
\]
so that a gerbe is a $1$-gerbe, and a bundle is a $0$-gerbe, in these terms. For $n\geq 2$, the $n$-gerbe itself is the $(0,n+2)$ part of $\cal H$, its multilayered connection consists of the $(1,n+1)$ to $(n+1,1)$ parts of $\cal H$, and the $n$-gerbe curvature is minus the globally-defined $(n+2)$-form part of $\cal H$. 

It is not particularly illuminating to write down the explicit equations for the higher gerbes. However it is amusing to consider a lower case, namely a 
$(-1)$-gerbe. This is, by the above definition,
\[
{\cal H} =\ln f_i - iA \in \Lambda ^{(1)}, \quad D{\cal H}=0,
\]
i.e. the $(-1)$-gerbe itself is a collection of functions $\ln f_i$, there is no layer of connections, and the curvature is the $1$-form $A$. Let us give an example of such a $(-1)$-gerbe on $S^1$. We cover the circle, parametrised by a $2\pi$-periodic coordinate $\theta$, with three open sets: $U_1=]0,\pi[$, $U_2=]2\pi/3,5\pi/3[$ and $U_3=]4\pi/3,7\pi/3[$. 
The equation $D {\cal H}=0$ implies the following equations:
\[
\ln f_j - \ln f_i = 0
\]
on $U_{ij}$,
\[
d\ln f_i -iA=0
\]
on $U_i$, and 
\[
dA=0
\]
on $M$, which are solved by 
\[
\ln f_i (\theta) = i\theta, \quad A=d\theta .
\]
In terms of the ``gerbopole'' example at the end of the previous section, the $(-1)$-gerbe sits on the equator of the monopole bundle, in the same way as the monopole bundle sits on the equator of the gerbopole.

The notion of gauge equivalence between bundles and gerbes in the previous two sections can also be expressed in terms of the language introduced in this section. For bundles the equations (\ref{be1}) and (\ref{be2}) for gauge equivalence can be expressed as:
\[
{\cal B'} \sim {\cal B} \Leftrightarrow {\cal B'}- {\cal B}=D \ln h_i, \quad
\ln h_i\in \Lambda^{0,1},
\]
and for gerbes the equations (\ref{ge1}) to (\ref{ge3}) for gauge equivalence can be expressed as:
\[
{\cal G'} \sim {\cal G} \Leftrightarrow {\cal G'}- {\cal G}=D (\ln h_{ij}+ i B_i), \quad
\ln h_{ij} + i B_i \in \Lambda^{0,2} \oplus \Lambda^{1,1}.
\]
This generalises to gauge equivalence for $n$-gerbes $\cal H'$, ${\cal H}\in \Lambda^{(n+2)}$: 
\[
{\cal H'} \sim {\cal H} \Leftrightarrow {\cal H'}- {\cal H}=D {\cal F}, \quad
{\cal F} \in \Lambda_0^{(n+1)},
\]
where $\Lambda_0^{(n+1)}$ denotes $\Lambda^{(n+1)}$ without the top-degree $(n+1)$-form part. 

We conclude this section by remarking that clearly it would be more natural if we could replace $\Lambda_0^{(n+1)}$ simply by $\Lambda^{(n+1)}$ in the above definition of equivalence. For example, in the case of bundles the higher gauge equivalence suggested here modifies equation (\ref{be2}) as follows, taking ${\cal F}= \ln h_i + i B$,
\[
i(A'_i-A_i) = d\ln h_i + iB,
\]
but also implies a higher gauge transformation for the curvature $F$:
\[
F'-F=dB.
\]
This would not normally be considered as a gauge transformation, but it is an equivalence for some purposes, since e.g. Chern forms for bundles on closed manifolds are preserved under this transformation.

\label{cohom}
\section{Comments}
Given the importance of abelian gauge theory, and the fact that gerbes generalise abelian bundles with connection in such a natural way, it would be very interesting if a direct, dynamical role could be found for gerbes in physical theories. For this it is necessary to construct couplings to other fields and build actions. A possible coupling arises in so-called twisted vector bundles, where abelian gerbes with connection are coupled to non-abelian bundles with connection (see \cite{m}, where the holonomy of such objects is also discussed). Actions for gerbes have been studied by Baez \cite{Bae}, in fact, in the context of non-abelian gerbes. Adapting the equations in section~\ref{locgerbe} to non-abelian gerbes is a challenging task --- see \cite{BreMes} for an algebraic/differential-geometric approach and \cite{Att} for a combinatorial approach. Maybe the cohomological approach and example presented here will suggest some way forward in this problem.

\label{Comm}


\bibliographystyle{plain}

\begin{thebibliography}{10}



\bibitem{Att}
{R. Attal.}
\newblock {Combinatorics of Non-Abelian Gerbes with Connection and Curvature}.
\newblock math-ph/0203056, 2002.



\bibitem{Bae}
{J. Baez}.
\newblock {Higher Yang-Mills Theory}.
\newblock hep-th/0206130, 2002.



\bibitem{Bar91}
J.~W. Barrett.
\newblock Holonomy and path structures in general relativity and {Y}ang-{M}ills
  theory.
\newblock {\em Int. J. Theor. Phys.}, 30(9):1171--1215, 1991.




\bibitem{BreMes}
{L. Breen and W. Messing}.
\newblock {Differential Geometry of Gerbes}.
\newblock math.AG/0106083, 2001.



\bibitem{Bry}
J-W. Brylinski.
\newblock {\em Loop spaces, characteristic classes and geometric quantization},
  volume 107 of {\em Progress in Mathematics}.
\newblock Birkhauser, 1993.


\bibitem{CP94}
A.~Caetano and R.~F. Picken.
\newblock An axiomatic definition of holonomy.
\newblock {\em Int. J. Math.}, 5(6):835--848, 1994.


\bibitem{cmm}
{A. Carey, J. Mickelsson and M. Murray.}
\newblock {Bundle gerbes applied to quantum field theory.}
\newblock {\em Rev. Math. Phys.}, 12(1):65--90, 2000.


\bibitem{Gir}
J.~Giraud.
\newblock {\em Cohomologie non-abelienne}, volume 179 of {\em Grundl.}
\newblock {S}pringer-{V}erlag, {B}erlin, 1971.

\bibitem{gom}
{K. Gomi.}
\newblock {Gerbes in classical Chern-Simons theory.}
\newblock hep-th/0105072, 2001.


\bibitem{Hi99}
{N.~Hitchin.}
\newblock Lectures on special {L}angrangian submanifolds.
\newblock In {\em {S}chool on {D}ifferential {G}eometry (1999)}, the {A}bdus
{S}alam {I}nternational {C}entre for {T}heoretical {P}hysics.

\bibitem{m}
{M. Mackaay.}
\newblock { A note on the holonomy of connections in twisted
bundles}.
\newblock To appear in {\em Cahiers Topologie G{\'e}om.
Diff{\'e}rentielle Cat{\'e}g.}, math.DG/0106019, 2001.



\bibitem{MP}
M. Mackaay and R. Picken.
\newblock {Holonomy and Parallel Transport for Abelian Gerbes}.
\newblock {\em {Adv. Math.}}, 170:287--339, 2002.


\bibitem{mei}
{E. Meinrenken}.
\newblock {The basic gerbe over a compact simple Lie group}.
\newblock math.DG/0209194, 2002.


\bibitem{P03}
{R. Picken.}
\newblock {TQFT's and gerbes.}
\newblock math.DG/0302065, 2003.


\bibitem{PS}
{R. Picken and P. Semi\~{a}o}.
\newblock {A classical approach to TQFT's}.
\newblock math.QA/0212310, 2002.


\bibitem{wit}
\newblock {E. Witten.}
\newblock { D-branes and $K$-theory.}
\newblock {\em J. High Energy Phys.}, no. 12, paper 19, 41 pp. (electronic), 1998.













\end{thebibliography}

\end{document}